\documentclass[11pt]{amsart}
\usepackage{amssymb,amsmath,amsfonts,amsthm,amscd}
\newtheorem{theorem}{Theorem}[section]
\newtheorem{lemma}[theorem]{Lemma}
\newtheorem{fact}[theorem]{Fact}
\newtheorem*{thm}{Theorem}
\newtheorem*{mthm}{Main Theorem}
\newtheorem*{mprop}{Theorem}
\newtheorem*{slemma}{Special Lemma}

\newtheorem*{definition}{Definition}

\newtheorem{prop}[theorem]{Proposition}
\newtheorem{cor}[theorem]{Corollary}

\newtheorem{remark}[theorem]{Remark}

\numberwithin{equation}{section}

\newcommand{\coker}{\mathop{\mathrm{coker}}\nolimits}
\newcommand{\codim}{\mathop{\mathrm{codim}}\nolimits}

\newcommand{\ann}{\mathop{\mathrm{Ann}}\nolimits}
\newcommand{\supp}{\mathop{\mathrm{Supp}}\nolimits}

\newcommand{\depth}{\mathop{\mathrm{depth}}\nolimits}
\newcommand{\ext}{\mathop{\mathrm{Ext}}\nolimits}

\newcommand{\onto}{\twoheadrightarrow}
\newcommand{\into}{\hookrightarrow}

\newcommand{\iso}{\simeq}
\newcommand{\im}{\mathop{\mathrm{Im}}\nolimits}

\newcommand{\projdim}{\mathop{\mathrm{proj\thinspace dim}}\nolimits}

\newcommand{\tor}{\mathop{\mathrm{Tor}}\nolimits}

\renewcommand{\l}{\ell}

\renewcommand{\hom}{\mathop{\mathrm{Hom}}\nolimits}

\renewcommand{\bar}{\protect\overline}
\renewcommand{\hat}{\protect\widehat}

\begin{document}

\title{THE UPPER BOUND OF FROBENIUS RELATED LENGTH FUNCTIONS}

\author{JINJIA LI}
\address{Department of Mathematics, University of Illinois at
Urbana-Champaign, Urbana, Illinois 61801}
\curraddr{Department of Mathematics, University of Illinois at
Urbana-Champaign, Urbana, Illinois 61801}
\email{jinjiali@math.uiuc.edu}
\date{September 10, 2004.}


\begin{abstract}
In this paper, we study the asymptotic behavior of lengths of
$\tor$ modules of homologies of complexes under the iterations of
the Frobenius functor in positive characteristic. We first give
upper bounds to this type of length functions in lower dimensional
cases and then construct a counterexample to the general
situation. The motivation of studying such length functions arose
initially from an asymptotic length criterion given in [D4] which
is a sufficient condition to a special case of nonnegativity of
$\chi_\infty$. We also provide an example to show that this
sufficient condition does not hold in general.
\end{abstract}

\maketitle

\specialsection*{INTRODUCTION AND NOTATIONS} In this paper, ($A$,
$m$, $k$) will be a complete local ring of characteristic $p>0$,
$m$ its maximal ideal, $k=A/m$ and $k$ is perfect. By a free
complex we mean a complex $F_\bullet = (F_i, d_i)_{i\geq 0}$
($\cdots\to F_2 \overset{d_2}{\to} F_1 \overset{d_1}{\to} F_0 \to
0$) of finitely generated free $A$-modules. We define codimension
of $M$ to be $\dim A - \dim M$ (denoted by $\codim M$) for any
$A$-module $M$. The Frobenius endomorphism $f_A:A \to A$ is
defined by $f_A(r)=r^p$ for $r \in A$. Each iteration $f_A^n$
defines a new $A$-module structure on $A$, denoted by
${}^{f^n}\!\! A$ for which $a\cdot b=a^{p^n}b$. Write $F_A^n(M)$
for $M \otimes_A {}^{f^n}\!\! A$ and $F_A^n(F_\bullet)$ for
$F_\bullet \otimes_A {}^{f^n}\!\! A$. We drop the subscript $A$
when there is no ambiguity.

In [D1], Dutta introduced the following definition of
$\chi_\infty$.

\begin{definition}
Let $R$ be a local ring in characteristic $p>0$. Let $M$ and $N$
be two finitely generated modules such that $\l(M \otimes_R
N)<\infty$ and $\projdim M < \infty$. Define
\[\chi_\infty(M,N)= \lim_{n\to \infty} \chi(F^n(M), N)/p^{n\codim
M}.\]
\end{definition}

$\chi_\infty$ plays an important role in the study of intersection
multiplicity $\chi$ defined by Serre ([S]), especially in the
non-smooth situation. For example, over complete intersections,
$\chi_\infty(M,N)=\chi(M,N)$ when both $M$ and $N$ are of finite
projective dimension ([D4], Corollary to Theorem 1.2). Thus the
positivity (or nonnegativity) of $\chi_\infty$ settles the
positivity (respectively, nonnegativity) conjecture of $\chi$ over
complete intersections.

Our main object is to examine the following sufficient condition
for the nonnegativity of $\chi_\infty$. ([D4], Corollary 1 to
Theorem 2.2)
\begin{thm}[Dutta] Let $R$ be a local Gorenstein ring in
characteristic $p>0$. Let $M$ and $N$ be finitely generated
modules of finite projective dimensions such that $\l(M \otimes N)
< \infty$. Suppose $\dim M + \dim N= \dim R$, $\dim N=\depth N +1
=s$ and $\dim M =\depth M +1 = 2$. Then $\chi_\infty(M,N) \geq 0$,
if
\[ \underset{n\to \infty}{\lim} \l(\ext^3(N,R)\otimes
H_m^0(F^n(\ext^{s+1}(M,R)))^\vee)/p^{ns}=0 . \tag{1}\]
\end{thm}
Note here $\projdim M=s+1$ and $\projdim N=3$ by
Auslander-Buchsbaum formula, and these $\ext$s are the natural
duals under the generalized ``Matlis" duality.

This study led us to investigate the asymptotic behavior of
$\l(\tor_j^A(H_i(F^n(F_\bullet)),N))$, where $F_\bullet$ is a free
complex with homologies of finite length. ($F_\bullet$ is not
necessarily a bounded complex here!).

In [D3], Dutta established that
\[\l(\tor_j^A(H_i(F^n(F_\bullet)),N))\leq
C_{ij}p^{n\dim N}\] when $\codim N =1$ ([D3], Proposition 1.3).
Naturally, one can ask whether this inequality is still valid when
$N$ has higher codimension. Investigation of the length condition
(1) raises the same question. The expectation was that the same
inequality should hold in general for any $N$, namely,
$\l(\tor_j^A(H_i(F^n(F_\bullet)),N))\leq C_{ij}p^{n\dim N}$. A
positive answer to this question in codimension 3 would yield an
affirmative answer for (1). However, our investigation revealed
that one can only extend this for $\codim N \leq 2$.

The following result in section 1 shows that one can extend this
inequality for $\codim N \leq 2$.
\begin{mprop}[Corollary 1.3 in section 1]
Let $F_\bullet$ be a free complex with homologies of finite length
over a Cohen-Macaulay local ring $A$. Let $N$ be a finitely
generated $A$-module such that $\codim N \leq 2$. Then there exist
constants $C_{ij}$'s, such that
\[ \l(\tor_j^A(H_i(F^n(F_\bullet)),N))\leq C_{ij}p^{n\dim N} \]
for all $i,j \geq 0$
\end{mprop}

When $\codim N=3$, we provide a counterexample in section 2. This
counterexample in turn leads to us our main theorem in section 3.
\begin{mthm}[Theorem 3.2 in section 3] Let $R=K[[X,Y,U,V]]/(XY-UV)$ where $K$ is
a field of characteristic $p>0$ and $X$, $Y$, $U$, $V$ are
indeterminates. There exist finitely generated modules $M$, $N$
over $R$ as in the above theorem with $s=1$, such that the
sufficient condition (1) for nonnegativity of $\chi_\infty$ fails
to hold.
\end{mthm}
Nevertheless, this counterexample does not give a negative
$\chi_\infty$.

\section{}
We first state a proposition due to Seibert ([Se], Proposition 1,
section 3) which plays a crucial role in our proof.
\begin{prop}[Seibert]Let $F_\bullet$ be a free complex over $A$ with homologies of finite length
and $N$ be any finitely generated $A$-module. Then there exist
constants $C_i$'s such that
\[ \l(H_i(F^n(F_\bullet)\otimes_A N))\leq C_{i}p^{n\dim N} \]
\end{prop}

The following is our first result which generalizes a result due
to Dutta ([D3], Proposition 1.3).
\begin{prop}
Let $F_\bullet$ be a free complex with homologies of finite length
over $A$. Let $N$ be $A/xA$ or $A/(x,y)$ where \{$x$\} or
\{$x$,$y$\}, respectively, forms a regular sequence. Then there
exist constants $C_{ij}$'s, such that
\[ \l(\tor_j^A(H_i(F^n(F_\bullet)),N))\leq C_{ij}p^{n\dim N} \tag{2}\]
for all $i,j \geq 0$.

\end{prop}

The following special lemma has been used repeatedly in the proof
of Proposition 1.2. We leave the proof as an exercise for the
reader.
\begin{slemma}Let $A$ be a local ring and $M$ be a module over $A$
such that $\l(M)<\infty$. Suppose $x$ is $A$-regular. Then
\[\l(\tor_1^A(M, A/xA)=\l(M \otimes_A (A/xA))\]
\end{slemma}

\begin{proof}[Proof of Proposition 1.2]
We write $\bar A=A/xA$ and $\bar F_\bullet =F_\bullet \otimes_A
\bar A$.\\
\textit{\underline {Case 1}}. $N=A/xA$.\\
This case has already been demonstrated in [D3] in a more general
set up. (See the proof of Proposition 1.3 in [D3], although the
official statement there is in the form of limit). We give a
simple proof of this case anyway for completeness.

Since $\projdim N=1$,
\[\tor_j^A(H_i(F^n(F_\bullet)),N)=0\] for $j\geq2$ and by the special lemma
\[\l(\tor_1^A(H_i(F^n(F_\bullet)),N))=\l(H_i(F^n(F_\bullet))\otimes_A
N).\] Thus it suffices to prove the result for $j=0$.\\
If $i=0$, since $H_0(F^n(F_\bullet))\otimes
N=H_0(F^n(F_\bullet)\otimes N)$, we get
the desired inequality by Proposition 1.1.\\
If $i\geq 1$, since $F_A^n(F_\bullet) \otimes_A \bar{A} =
F_{\bar{A}}^n(\bar{F_\bullet})$, there is a short exact sequence
of complexes
\[ 0
        \to F^n(F_\bullet)
        \overset{x}{\to} F^n(F_\bullet)
        \to F_{\bar{A}}^n(\bar{F_\bullet})
        \to 0 .\]
Taking the associated long exact sequence of homologies, we get
\[\cdots
    \to H_i(F^n(F_\bullet))
    \overset{x}{\to} H_i(F^n(F_\bullet))
    \to H_i(F_{\bar{A}}^n(\bar{F_\bullet}))
    \to H_{i-1}(F^n(F_\bullet))
    \to \cdots .\]
It yields the following short exact sequence
\[ 0
      \to H_i(F^n(F_\bullet)) \otimes A/xA
      \to H_i(F_{\bar{A}}^n(\bar{F_\bullet}))
      \to (0:x)_{H_{i-1}(F^n(F_\bullet))}
      \to 0 \tag{3}
\]
for $i \geq 1$. So,
\[\l(H_i(F^n(F_\bullet)) \otimes A/xA)\leq
\l(H_i(F_{\bar{A}}^n(\bar{F_\bullet})))\] and again, the desired
inequality follows from Proposition 1.1 with $N=\bar{A}$.

\textit{\underline{Case 2}}. $N=A/(x,y)$.\\
In this case, since $\projdim N=2$,
\[\tor_j^A(H_i(F^n(F_\bullet)),N)=0\] for $j\geq3$. By a result
due to Serre ([S], Theorem 1, Chapter IV),
\[ \sum_{j=0}^{2}(-1)^j \l(\tor_j^A(H_i(F^n(F_\bullet)),N))=\chi(H_i(F^n(F_\bullet)),N))=0. \]
Hence, it is enough to prove the result for $j=0$ and $1$.\\
Tensor (3) with $A/(x,y)$ ($\iso \bar A/y \bar A$) over $\bar A$.
We obtain
\[ \cdots
           \to \tor_1^{\bar A}((0:x)_{H_{i-1}(F^n(F_\bullet))},{\bar A}/y{\bar A})
           \to H_i(F^n(F_\bullet)) \otimes_A A/(x,y) \to \]
           \[\to H_i(F_{\bar{A}}^n(\bar{F_\bullet})) \otimes_{\bar A} {\bar A}/y{\bar A}
           \to  (0:x)_{H_{i-1}(F^n(F_\bullet))} \otimes_{\bar A}{\bar A}/y{\bar A}
           \to 0  \tag{4}    \]
for $i \geq 1$. It follows that
\begin{align*}
    \l(H_i(&F^n(F_\bullet)) \otimes_A A/(x,y))\\
                 \leq& \l(\tor_1^{\bar A}((0:x)_{H_{i-1}(F^n(F_\bullet))},{\bar A}/y{\bar A}))
                   +\l(H_i(F_{\bar{A}}^n(\bar{F_\bullet})) \otimes_{\bar A} {\bar A}/y{\bar A}).
\end{align*}
Notice that by the special lemma
\[\l(\tor_1^{\bar A}((0:x)_{H_{i-1}(F^n(F_\bullet))},{\bar A}/y{\bar A}))
=\l((0:x)_{H_{i-1}(F^n(F_\bullet))} \otimes_{\bar A}{\bar
A}/y{\bar A})\] and from the above long exact sequence (4)
\[\l((0:x)_{H_{i-1}(F^n(F_\bullet))} \otimes_{\bar A}{\bar A}/y{\bar
A}) \leq \l(H_i(F_{\bar{A}}^n(\bar{F_\bullet})) \otimes_{\bar A}
{\bar A}/y{\bar A}).\] Hence
\[\l(H_i(F^n(F_\bullet)) \otimes_A A/(x,y))
                 \leq 2\l(H_i(F_{\bar{A}}^n(\bar{F_\bullet})) \otimes_{\bar A} {\bar A}/y{\bar
                 A})\]
Therefore by case 1, we are done for $j=0$.\\
Finally, for $j=1$, we use the following spectral sequence
obtained by base change,
\[  \tor_p^{\bar A}(\tor_q^A(H_i(F^n(F_\bullet)) ,
           \bar A), {\bar A}/y{\bar A})
    \underset {p} {\Longrightarrow}
    \tor_{p+q}^A(H_i(F^n(F_\bullet)) , A/(x,y)).
\]
It follows that
\begin{align*}
\l(\tor_1^A&(H_i(F^n(F_\bullet)) , A/(x,y)))\\
\leq& \l(\tor_1^{\bar A}(\tor_0^A(H_i(F^n(F_\bullet)) ,
           \bar A), {\bar A}/y{\bar A}))\\
    &\ \ \ \ \ +    \l(\tor_0^{\bar A}(\tor_1^A(H_i(F^n(F_\bullet)) ,
           \bar A), {\bar A}/y{\bar A}))\\
=&  \l(H_i(F^n(F_\bullet)) \otimes A/(x,y))\\
    &\ \ \ \ \ +\l(\tor_1^A(H_i(F^n(F_\bullet)), \bar A) \otimes
    A/(x,y)).
\end{align*}
The last equality here is by the special lemma again.\\
Since $x$ is $A$-regular, $\tor_1^A(H_i(F^n(F_\bullet)), \bar A)
\iso (0:x)_{H_i(F^n(F_\bullet))}$. Therefore by (3), we have a
surjection
\[H_{i+1}(F_{\bar{A}}^n(\bar{F_\bullet}))\otimes A/(x,y) \onto
    \tor_1^A(H_i(F^n(F_\bullet)), \bar A) \otimes A/(x,y).
\]
Thus,
\begin{align*}
  \l(\tor_1^A&(H_i(F^n(F_\bullet)) , A/(x,y)))\\
    \leq& \l(H_i(F^n(F_\bullet)) \otimes A/(x,y))
    +\l(H_{i+1}(F_{\bar{A}}^n(\bar{F_\bullet}))\otimes
        A/(x,y)).\\
\end{align*}
Both of the terms on the right hand side of the above inequality
are bounded by a constant times $p^{n\dim N}$ by the $j=0$ case,
and so we are done for $j=1$ which finishes our proof.
\end{proof}

\begin{cor}Let $A$ be a Cohen-Macaulay local ring and let $F_\bullet$ be as in Proposition 1.2.
Let $N$ be a finitely generated $A$-module such that $\codim N
\leq 2$. Then there exist constants $C_{ij}$'s, such that
\[ \l(\tor_j^A(H_i(F^n(F_\bullet)),N))\leq C_{ij}p^{n\dim N} \]
for all $i,j \geq 0$
\end{cor}
\begin{proof}
Suppose $\codim N=h$, $h=1$ or $2$. Then $\ann_A N$ contains an
$A$-regular sequence \{$x_1$,...,$x_h$\}. We have the following
short exact sequence,
\[0\to Q\to (A/(x_1,...,x_h))^t \to N \to 0\]
Tensoring the above short exact sequence with
$H_i(F^n(F_\bullet))$, we get a long exact sequence
\[\cdots \to \tor_1^A((A/(x_1,...,x_h))^t,H_i(F^n(F_\bullet))) \to
\tor_1^A(N,H_i(F^n(F_\bullet))) \to Q \otimes H_i(F^n(F_\bullet))
\to \]
\[ \to (A/(x_1,...,x_h))^t \otimes H_i(F^n(F_\bullet)) \to N
\otimes H_i(F^n(F_\bullet)) \to 0 .\] By Proposition 1.2 and
induction on $j$, we obtain the desired inequality.
\end{proof}
\begin{remark} If $A$ is a regular local ring, since the functor
$F^n(-)$ is exact ([K], Theorem 3.3),
$\tor_j^A(H_i(F^n(F_\bullet)),N) \iso
\tor_j^A(F^n(H_i(F_\bullet)),N)$. Thus by Proposition 1.1 the
inequality in Proposition 1.2 hold for any finitely generated
$A$-module $N$ .
\end{remark}

\section{}
Now, we demonstrate an example to show that the inequality (2) in
Proposition 1.2, as well as the one in Corollary 1.3, can fail
when $\codim N=3$.

We first state two standard facts in commutative algebra which
will be used in the proof of Proposition 2.4.
\begin{fact}
Let $R$ be a finitely generated algebra over a field $K$ and $M$
be a finitely generated $R$-module. Let $m$ be a maximal ideal of
$R$. Suppose $\supp{M}=\{ m \}$ and $K \iso R/m$ via the natural
map. Then $\l_R(M)=\dim_K M$. Here $\dim_K M$ denote the dimension
of $M$ as a $K$-vector space.
\end{fact}
\begin{fact}
Let $R$ be a commutative ring and $M$ be a finitely generated
$R$-module. Let $\hat{R_m}$ be the $m$-adic completion of $R_m$
where $m$ is a maximal ideal of $R$. If\ \ $\supp_R M=\{m\}$, then
\[\l_R (M)= \l_{R_m} (M_m)=\l_{\hat{R_m}} (\hat{M_m}).\]
\end{fact}

\begin{lemma}
Let $R=K[X,Y,U,V]/(XY-UV)$ where $K$ is a field of characteristic
$p>0$ and $X$, $Y$, $U$, $V$ are indeterminates. Consider $K$ as a
module over $R$ in the obvious way. Then $\hom_R(K, F_R^n(K))$ is
a $K$-vector space and
\[\dim_K \hom_R(K, F_R^n(K)) \geq p^n \]
\end{lemma}

\begin{proof}
To simplify our notations, we use $x$,$y$,$u$,$v$ to denote the
images of $X$, $Y$, $U$, $V$ respectively in any quotient ring of
$K[X,Y,U,V]$ if there is no confusion about that ambient quotient
ring. $\hom_R(K, F_R^n(K))$ is a $K$-vector space consisting of
all the elements of $F_R^n(K)$ which are killed by the maximal
ideal $(x,y,u,v)$. Let $\mathcal A = \{x^{p^n-1}y^iu^{p^n-1-i}| 0
\leq i \leq p^n-1\}$, which is a subset of
\[F_R^n(K)=\frac{K[X,Y,U,V]}{(X^{{}^{p^n}}, Y^{{}^{p^n}},
U^{{}^{p^n}}, V^{{}^{p^n}},XY-UV)}.\] It is easy to verify that
$\mathcal A \subset \hom_R(K, F_R^n(K))$. We will show that
elements in $\mathcal A$ are
linearly independent over $K$ which gives us the desired inequality.\\

Let $\{\lambda_i\}_{0\leq i \leq p^n-1}$ be elements in $K$ such
that
\[ \sum_{i=0}^{p^n-1} \lambda_ix^{p^n-1}y^iu^{p^n-1-i}=0\in F_R^n(K) \tag{5}.\]
Let
\[S=\frac{K[X,Y,U,V]}{(X^{{}^{p^n}}, Y^{{}^{p^n}},
U^{{}^{p^n}}, V^{{}^{p^n}})}.\] Then $R=S/(xy-uv)$. Lift the
relation (5) to a relation in $S$. Since $S$ is a $K$-vector space
with basis \{$x^iy^ju^kv^l|0\leq i,j,k,l\leq p^n-1$ \}, we obtain
\[ \sum_{i=0}^{p^n-1}
\lambda_ix^{p^n-1}y^iu^{p^n-1-i}=(\sum_{0\leq i,j,k,l\leq
p^n-1}\mu_{{}_{i,j,k,l}} x^iy^ju^kv^l)(xy-uv)\in S \tag{6},\]
where the $\mu_{{}_{i,j,k,l}}$ are elements of $K$. Define
\[\lambda_{{}_{i,j,k,l}}=
    \begin{cases}
        \lambda_j,& \text{if $i=p^n-1$, $j+k=p^n-1$ and $l=0$},\\
        0,& \text{otherwise}.
    \end{cases}
\]
We also define $\mu_{{}_{i,j,k,l}}=0$ if one of
$i$, $j$, $k$, $l$ is negative. \\
By comparing the coefficients on both sides of $(6)$, we obtain
that
\[\lambda_{{}_{i,j,k,l}}=\mu_{{}_{i-1,j-1,k,l}}-\mu_{{}_{i,j,k-1,l-1}}, \ \ \
\forall \ \ i,j,k,l \leq p^n-1.
\]
Using the above formula repeatedly, noticing that
$\lambda_{{}_{i,j,k,l}}=0$ if $i<p^n-1$, we get
\begin{align*}
\lambda_i &= \lambda_{p^n-1,i,p^n-1-i,0}\\
          &=\mu_{p^n-2,i-1,p^n-1-i,0}\ \ +\ \ 0 \\
          &=\mu_{p^n-3,i-2,p^n-i,1} \\
          &=\mu_{p^n-4,i-3,p^n-i+1,2} \\
          &\ \vdots\\
          &=\mu_{p^n-i-1,0,p^n-2,i-1} \\
          &=0 
\end{align*}
for all $i=0$, $1$, ..., $p^n-1$.
\end{proof}

The following is an example where the inequality (2) in
Proposition 1.2 fails when $\codim N =3$. The complex $F_\bullet$
is taken to be a free resolution of $K$ and $i=0$.
\begin{prop}
Let $R=K[[X,Y,U,V]]/(XY-UV)$ where $K$ is a field of
characteristic $p>0$ and $X$, $Y$, $U$, $V$ are indeterminates.
Then
\[\l(\tor_3^R(F^n(K), R/(x,y,u+v))) \geq p^n.\]
\end{prop}

\begin{proof}
Since \{$x$, $y$, $u+v$\} forms an $R$-sequence, it follows that
\[\tor_3^R(F^n(K), R/(x,y,u+v)) \iso \hom_R( R/(x,y,u+v),F^n(K)). \]
Since there is a surjection $R/(x,y,u+v)\onto K$, by applying
$\hom_R(-,F^n(K))$, we obtain a injection
\[\hom_R(K, F^n(K))\into \hom_R( R/(x,y,u+v),F^n(K)).\]
From Fact 2.1, Fact 2.2 and Lemma 2.3, we have
\[\l(\hom_R(K, F^n(K))) \geq p^n.\]
Therefore
\[\l(\tor_3^R(F^n(K), R/(x,y,u+v))) \geq p^n.\]
\end{proof}

\begin{remark}
Using the same method one can show that over the hypersurface ring
$R=K[[X_1,...,X_t, Y_1,...,Y_t]]/(\sum_{i=1}^t X_iY_i)$,
$\l(\hom_R(K, F^n(K)))$ is unbounded.
\end{remark}

\section{}
In [D4], Dutta gave an asymptotic length condition over Gorenstein
local rings of positive characteristic for the nonnegativity of
$\chi_\infty(M,N)$ when $\dim M = 2$. In this section, we will
construct examples to show that over the local hypersurface $R$
discussed in Corollary 2.4, this length condition fails to hold.

Let $R$ be a local ring in characteristic $p > 0$. Let $M$ and $N$
be two finitely generated modules such that $\l(M\otimes_R N) <
\infty$, $\dim M + \dim N \leq \dim R$ and $\projdim M < \infty$.
In [D1], Dutta defined
\[ \chi_\infty(M,N)=\underset{n\to
\infty}{\lim}\chi(F^n(M),N)/p^{n\codim M} .\] For properties of
$\chi_\infty$, see [D1, D2, R, Se]. Dutta [D4] established the
following criterion for nonnegativity of $\chi_\infty$ over a
local Gorenstein rings of positive characteristic.

\begin{theorem}[Dutta]Let $R$ be a local Gorenstein ring in
characteristic $p>0$. Let $M$ and $N$ be finitely generated
modules of finite projective dimension such that $\l(M \otimes N)
< \infty$. Suppose $\dim M + \dim N= \dim R$, $\dim N=\depth N +1
=s$ and $\dim M =\depth M +1 = 2$. Then $\chi_\infty(M,N) \geq 0$,
if
\[ \underset{n\to \infty}{\lim} \l(\ext^3(N,R)\otimes
H_m^0(F^n(\ext^{s+1}(M,R)))^\vee)/p^{ns}=0 .\]

Here, $(-)^\vee $ denotes the Matlis duality $\hom_R(-,E)$ where
$E$ is the injective hull of the residue field of $R$.
\end{theorem}

The following is an example where the length criterion in Theorem
3.1 fails.

\begin{theorem}
Let $R=K[[X,Y,U,V]]/(XY-UV)$ where $K$ is a field of
characteristic $p>0$ and $X$, $Y$, $U$, $V$ are indeterminates.
There exist finitely generated modules $M$, $N$ over $R$
satisfying the conditions in Theorem 3.1 such that
\[ \underset{n\to \infty}{\lim} \l(\ext^3(N,R)\otimes
H_m^0(F^n(\ext^2(M,R)))^\vee)/p^n>0 .\]
\end{theorem}

\begin{proof}
We are going to construct modules $M$ and $N$ satisfying the
conditions in Theorem 3.1 with $s=1$, such that $\ext_R^2(M,R)
\iso K$ and $\ext_R^3(N,R) \iso K$.

Let $x$, $y$, $u$, $v$ denote the images of $X$, $Y$, $U$, $V$ in
$R$. Take a minimal free resolution of $K$ over $R$
\[ \cdots \to R^t \overset{\psi}{\to} R^4 \overset{\phi}{\to} R \to K \to 0 \]
where $\phi$ can be written as a matrix $[x,y,u,v]$ with respect
to the standard bases for $R^4$ and $R$. Let $(-)^*$ denote
$\hom_R(-,R)$. Apply $(-)^*$ to the above exact sequence. Since
$\depth R =3$, $K^*=0$ and we obtain the following exact sequence
\[ 0 \to R \overset {\phi^*}{\to} R^4 \to M' \to 0\]
where $M'= \coker {\phi^*}$. Let \{$\mathbf{e}_1$, $\mathbf{e}_2$,
$\mathbf{e}_3$, $\mathbf{e}_4$\} be a standard basis for $R^4$, it
follows that
$M'=R^4/R(x\mathbf{e}_1+y\mathbf{e}_2+u\mathbf{e}_3+v\mathbf{e}_4)$.
Note that if $r \in \ann_R M'$, then there exists an $a \in R$
such that
\[r(\mathbf{e}_1+\mathbf{e}_2+\mathbf{e}_3+\mathbf{e}_4)=a(x\mathbf{e}_1+y\mathbf{e}_2+u\mathbf{e}_3+v\mathbf{e}_4).\]
It follows that $ax=ay=r$. But $R$ is a domain and $x \ne y$ in
$R$, thus $r=0$. Therefore $\ann_R M'=(0)$ whence $\dim M'= \dim R
=3$. Moreover, since $\ext_R^1(K, R)=0$, $M'=\im {\psi^* }$, which
is a submodule of $R^t$ and therefore torsion-free. Hence, $x \in
R$ is a non zero divisor on $M'$.

Let
$M=M'/xM'=R^4/(xR^4+R(x\mathbf{e}_1+y\mathbf{e}_2+u\mathbf{e}_3+v\mathbf{e}_4))$.
It follows that $\dim M=2$. One can also prove that $\projdim M=2$
since $\projdim M'=1$ and $x$ is both $M'$-regular and $R$-regular.
Therefore by Auslander-Buchsbaum formula, $\depth M=1$. Moreover,
\[\ext_R^2(M,R)\iso \ext_R^1(M',R) \iso K.\]

In order to construct $N$, let $\bar R= R/(y, u+v)$. Then $\dim
\bar R=1$, $\depth \bar R =1$. Take a minimal resolution of $K$
over $\bar R$
\[ {\bar R}^2 \overset{\zeta}{\to} \bar R \to K \to 0 .\]
Apply $\hom_{\bar R}(-, \bar R)$. Let $N=\coker{\zeta}^*$ and we
obtain a free resolution of $N$ over $\bar R$
\[ 0 \to \bar R \overset{\zeta^*}{\to} {\bar R}^2 \to N \to 0.
\]
Use a similar argument as before, $\ann_{\bar R}N =(0)$. Hence
$\dim_{\bar R} N=1$, $\projdim_{\bar R} N =1$ and $\depth_{\bar R}
N=0$. Therefore $\dim_R N=1$, $\depth_R N=0$ and $\projdim_R N=3$.
Note that $\l(M \otimes_R N) < \infty$ since the annihilator of $M
\otimes_R N$ contains $(x,y,u+v)$ which is primary to the maximal
ideal $(x,y,u,v)$. Moreover, $\ext_R^3(N,R) \iso \ext_{\bar
R}^1(N, \bar R) \iso K$.

Finally, to check
\[ \underset{n\to \infty}{\lim} \l(\ext^3(N,R)\otimes
H_m^0(F^n(\ext^2(M,R)))^\vee)/p^n>0 \]
It is enough to notice that
\begin{align*}
    &\l(\ext^3(N,R)\otimes H_m^0(F^n(\ext^2(M,R)))^\vee)\\
    =&\l((\ext^3(N,R)\otimes H_m^0(F^n(\ext^2(M,R)))^\vee)^\vee)\\
    =&\l(\hom(\ext^3(N,R),H_m^0(F^n(\ext^2(M,R)))))\\
    =&\l(\hom(\ext^3(N,R),H_m^0(F^n(K))))\\
    =&\l(\hom(K,F^n(K)))\\
    \geq&p^n\\
\end{align*}
\end{proof}

\begin{remark}
Although the length criterion does not hold in general, there do
exist local Gorenstein rings such that the length criterion holds
for all $M$ and $N$. It would be nice to have a general method to
identify such rings.
\end{remark}

\medskip

\specialsection*{ACKNOWLEDGEMENT} I am indebted to my thesis
advisor Sankar Dutta for his direction and many inspiring
discussions on the subject of this paper. I also would like to
thank the referee for the valuable suggestions and comments.

\bibliographystyle{amsalpha}

\end{document}